\documentclass[12pt,reqno]{amsart}

\usepackage{amsmath,amssymb,mathtools}
\usepackage{enumerate}
\usepackage{microtype}
\usepackage[usenames]{color} 
\usepackage[a4paper,margin=1.15in]{geometry}
\allowdisplaybreaks
\numberwithin{equation}{section}

\newtheorem{theorem}{Theorem}[section]
\newtheorem{Prop}[theorem]{Proposition}
\newtheorem{Lem}[theorem]{Lemma}
\theoremstyle{definition}
\newtheorem{Def}[theorem]{Definition}
\theoremstyle{remark}
\newtheorem{Rem}[theorem]{Remark}

\newcommand{\R}{\mathbb R}
\newcommand{\N}{\mathbb N}
\newcommand{\supp}{\operatorname{supp}}

\begin{document}
\raggedbottom

\title[Spectrality of alternating-sign Moran measures]{Spectrality of alternating-sign Moran measures}

\author{Jun Jason Luo}
\address{College of Mathematics and Statistics, Key Laboratory of Nonlinear
Analysis and its Applications (Chongqing University), Ministry of Education,
Chongqing University, Chongqing 401331, P.R. China}
\email{jun.luo@cqu.edu.cn}

\author{Lin Mao}
\address{College of Mathematics and Statistics, Key Laboratory of Nonlinear
Analysis and its Applications (Chongqing University), Ministry of Education,
Chongqing University, Chongqing 401331, P.R. China}
\email{maol@stu.cqu.edu.cn}

\author{Jing-cheng Liu}
\address{Key Laboratory of Computing and Stochastic Mathematics (Ministry of
Education), School of Mathematics and Statistics, Hunan Normal University,
Changsha, Hunan 410081, P.R. China}
\email{jcliu@hunnu.edu.cn}

\begin{abstract}
	For \(m\geq 2\), let \(D_m=\{0,1,\ldots,m-1\}\). For each \(k\geq 1\),
	consider the family of contractions
	\[
	\Phi_k=\{\phi_{k,d}:d\in D_{n_k}\},
	\qquad
	\phi_{k,d}(x)=(-1)^d b_k^{-1}(x+d),
	\]
	where \(b_k\) and \(n_k\) are integers satisfying \(b_k\geq n_k\geq 2\).
	We construct the canonical pullback attractor generated by
	\(\{\Phi_k\}_{k\geq 1}\) and the associated equal-weight Moran measure
	\(\mu\). Suppose that \(\{b_k\}_{k\geq 1}\) is bounded and that \(n_k\) is
	even for every \(k\). We prove that \(\mu\) is spectral if \(2\mid b_2,\ n_2\mid 2b_2$, and \(n_k\mid b_k\) for all	\(k\geq 3\).
	If, in addition, \(4\mid n_k\) for every \(k\), then the converse also holds,
	yielding a complete characterization of spectrality. The proof reduces the
	associated matrix-valued Fourier recursion to the Fourier transform of an
	infinite convolution of discrete probability measures. The main new ingredient
	in the necessity argument is a decomposition method for rational-scale infinite
	convolutions, in which a spectrum is successively partitioned into congruence
	classes.
\end{abstract}


\thanks{The research was supported by the National Natural Science Foundation of China
(No.~12071125), the Natural Science Foundation of Chongqing
(No.~CSTB2023NSCQ-MSX0553), the Hunan Provincial Natural Science Foundation
(No.~2024JJ3023), and the Scientific Research Fund of the Hunan Provincial
Education Department (No.~23A0059).}



\maketitle

\section{Introduction}
A compactly supported Borel probability measure \(\mu\) on \(\R^d\) is called
\emph{spectral} if there exists a countable set \(\Lambda\subset\R^d\) such
that
\[
E(\Lambda)
:=
\bigl\{e^{2\pi i\langle\lambda,\,\cdot\rangle}:\lambda\in\Lambda\bigr\}
\]
forms an orthonormal basis for \(L^2(\mu)\). In this case, \(\Lambda\) is
called a \emph{spectrum} of \(\mu\). When \(\mu\) is normalized Lebesgue
measure on an interval, this definition recovers the classical Fourier basis.
More generally, for  Lebesgue measure on a measurable set of finite
positive measure, determining spectrality is precisely the spectral-set
problem introduced by Fuglede \cite{Fuglede1974}. The notion of a spectral
measure thus extends the classical Fourier basis from Lebesgue measure to
certain singular measures, including measures supported on fractal sets.

Spectrality comprises two separate requirements. Orthogonality is governed by
the zero set of the Fourier transform: \((\Lambda-\Lambda)\setminus\{0\}
\subset \mathcal Z(\widehat\mu)\), whereas completeness requires \(E(\Lambda)\) to span \(L^2(\mu)\). For
singular measures, neither requirement can in general be determined from the
geometry of the support alone. The distribution of mass and the arithmetic
compatibility between contraction factors and digit sets also play essential
roles. The spectral analysis of fractal measures therefore involves a
combination of harmonic analysis, fractal geometry, and the arithmetic
structure of iterated function systems.

Jorgensen and Pedersen \cite{Jorgensen-Pedersen1998} constructed the first singular non-atomic spectral
measure, the middle-fourth Cantor measure. Their work initiated the systematic study of
orthogonal exponential bases for self-similar measures. For Cantor-type measures, finite-level orthogonality is expressed by the
unitarity condition defining a Hadamard triple. This principle led to the
constructions in \cite{Laba-Wang2002} and to the theorem that self-affine
measures associated with Hadamard triples are spectral
\cite{Dutkay-Han-Lai2019}.  In one dimension, the spectrality of Bernoulli convolutions was characterized by Hu and Lau \cite{Hu-Lau2008} and by Dai  \cite{Dai2012}, and Dai, He, and Lau \cite{Dai-He-Lau2014} obtained the corresponding
classification for self-similar measures generated by consecutive digit sets.

Moran constructions generalize the self-similar setting by allowing the
contraction ratios and digit sets to vary from level to level. The associated
measures are naturally represented as infinite convolutions of discrete
measures and, in general, do not satisfy a fixed self-similar identity. An early study of Fourier expansions associated with such measures was given by Strichartz \cite{Strichartz2000}. An and He \cite{An-He2014} established a divisibility condition
sufficient for the spectrality of Moran measures with consecutive digit sets. Subsequent work established sufficient conditions for broader classes of Cantor--Moran measures and infinite convolutions \cite{An-He-Lau2015,He-He2017,Shi2019}, and obtained necessary conditions or complete characterizations in several structured settings \cite{An-He-Li2015,Deng2022,Deng2023,An-Li-Zhang2025,Wu-Xiao2024}.
Meanwhile, tail measures, integral periodic zero sets, and weak convergence
were developed as general tools for proving spectrality of infinite
convolutions \cite{An-Fu-Lai2019,Li-Miao-Wang2022, Li-Miao-Wang2024-1,Li-Miao-Wang2024-2}.

Most results on spectral Moran measures allow the contraction ratio and the
digit set to vary with the level, but require all maps at a given level to
have the same linear part. A separate line of research concerns self-similar
measures with contractions of different signs. Wu and Liu \cite{Wu-Liu2022} constructed
examples of spectral self-similar measures with both positive and negative
contraction ratios. Wu \cite{Wu2024} later obtained a sharp
spectrality criterion for a self-similar system with a consecutive digit set
of even cardinality and alternating signs. More recently,
periodic sign patterns have been studied in the self-similar setting
\cite{Liu-Wang2025}. In these works, however, the system remains unchanged
from level to level.

The present paper combines the level dependence of Moran constructions with
alternating signs among the branches. Since the family of maps changes with
the level, the attractor is not obtained as the fixed point of a single
Hutchinson operator. We instead construct a canonical pullback attractor and
its associated equal-weight probability measure. Moreover, the alternating
signs prevent a direct use of the usual product formula for the Fourier
transform and lead to a matrix-valued recursion. By exploiting its rank-one
structure, we convert this recursion into the Fourier transform of an infinite
convolution of discrete probability measures.

Throughout the paper, we write \(D_m:=\{0,1,\ldots,m-1\}\) for \(m\geq1\).
For a Borel map \(S\), the pushforward of a measure \(\eta\) is denoted by
\(S_\#\eta\).

\begin{Def}\label{def-moran IFS}
Let $\{b_k\}_{k\geq1}$ and $\{n_k\}_{k\geq1}$ be sequences of integers with
$b_k\geq n_k\geq2$. For each $k\geq1$, let
$\Phi_k=\{\phi_{k,i}\}_{i\in D_{n_k}}$, where
\begin{equation}\label{ifs}
 \phi_{k,i}(x)=(-1)^i b_k^{-1}(x+i),\qquad x\in\R,\quad
 i\in D_{n_k}.
\end{equation}
We call $\{\Phi_k\}_{k\geq1}$ a Moran iterated function system with
alternating-sign linear parts.
\end{Def}

\begin{theorem}\label{thm-main1}
Let $\{\Phi_k\}_{k\geq1}$ be as in Definition \ref{def-moran IFS}.
\begin{enumerate}[(i)]
\item There exists a unique uniformly bounded sequence of nonempty compact
sets $\{E_k\}_{k\geq1}$ satisfying
\begin{equation}\label{equi-supp}
 E_k=\Phi_k(E_{k+1}):=\bigcup_{i=0}^{n_k-1}\phi_{k,i}(E_{k+1}),\qquad k\geq1.
\end{equation}
Moreover,
\begin{equation}\label{conv-point}
\begin{split}
&\lim_{l\to\infty}\bigl|\phi_{k,j_k}\circ\cdots\circ\phi_{l,j_l}(0)
-\phi_{k,j_k}\circ\cdots\circ\phi_{l,j_l}(a_{l+1})\bigr|=0
\end{split}
\end{equation}
uniformly in the admissible digits $j_i$ and in $a_{l+1}\in E_{l+1}$.

\item There exists a unique sequence of Borel probability measures
$\{\mu_k\}_{k\geq1}$ such that $\supp\mu_k=E_k$ and
\begin{equation}\label{equi-mea}
 \mu_k=\frac1{n_k}\sum_{i=0}^{n_k-1}(\phi_{k,i})_\#\mu_{k+1},
 \qquad k\geq1.
\end{equation}
\end{enumerate}
\end{theorem}

We call $\mu:=\mu_1$ the Moran measure associated with the system. Our two
spectrality results are as follows.

\begin{theorem}\label{thm-main2}
Assume that $\{b_k\}_{k\geq1}$ is bounded and that $2\mid n_k$ for every
$k\geq1$. Then $\mu$ is spectral provided that \(2\mid b_2, \  n_2\mid2b_2, \  n_k\mid b_k\) for \(k\geq3\).
\end{theorem}

Writing $n_k=2p_k$, the second condition above is equivalent to
$p_2\mid b_2$. In general these sufficient conditions are not necessary; see
Theorem \ref{Bernoulli situation}. When every $p_k$ is even, however, they are
also necessary.

\begin{theorem}\label{thm-main3}
Assume that $\{b_k\}_{k\geq1}$ is bounded and that $4\mid n_k$ for every
$k\geq1$. Then $\mu$ is spectral if and only if $ n_2\mid2b_2$ and $n_k\mid b_k$ for $k\geq 3$.
\end{theorem}

No arithmetic condition is required for \(b_1\), since changing \(b_1\)
merely dilates the measure, and spectrality is invariant under dilation.
The proof has two main steps. First, because the orientation in \eqref{ifs}
depends on the digit, the measure does not directly admit the usual
infinite-convolution representation. Instead, the real and imaginary parts of
its Fourier transform satisfy a recursion involving \(2\times2\) matrices.
These matrices have rank one, which converts the recursion into an ordinary
product formula and leads to an explicit infinite convolution that is spectral
if and only if \(\mu\) is spectral; see Proposition \ref{equivalent}.
Second, the necessity proof requires spectrality to be passed successively to
tail convolutions whose scale factors are rational rather than necessarily
integral. For this purpose, we develop a method for rational-scale infinite
convolutions that passes from a spectrum of the full convolution to spectra
of its successive tails by decomposing the spectrum into congruence classes.
Applying this method level by level yields the required divisibility
conditions.

The paper is organized as follows. Section \ref{sect.2} collects the spectral
criteria and infinite-convolution results used later. Section \ref{sect.3}
constructs the pullback attractor and proves Theorem \ref{thm-main1}. Section
\ref{sect.4} gives the Fourier reduction and proves the sufficient condition.
Section \ref{sect.5} develops the rational-scale decomposition and proves the
necessity statement.

\section{Preliminaries}\label{sect.2}

Let \(\mu\) be a Borel probability measure on \(\R\). We use the Fourier
transform convention
\[
 \widehat\mu(t)=\int_{\R}e^{2\pi itx}\,d\mu(x),\qquad t\in\R,
\]
and write
\(\mathcal Z(\widehat\mu)=\{t\in\R:\widehat\mu(t)=0\}\).
For a nonempty finite set \(A\subset\R\), we use the notation
\[
 \delta_A=\frac1{\#A}\sum_{a\in A}\delta_a,
 \qquad
 m_A(t)=\widehat{\delta_A}(t)
 =\frac1{\#A}\sum_{a\in A}e^{2\pi iat}.
\]
For finite sets \(A,B\subset\R\), we write
\(A+B=\{a+b:a\in A,\ b\in B\}\); the notation \(A\oplus B\) means that
every element of \(A+B\) has a unique such representation.
For a countable set \(\Lambda\subset\R\), let \(E(\Lambda)=\{e^{2\pi i\lambda x}:\lambda\in\Lambda\}\).
This family is orthogonal in \(L^2(\mu)\) if and only if \( (\Lambda-\Lambda)\setminus\{0\}\subset\mathcal Z(\widehat\mu)\).
Such a set \(\Lambda\) is called a \emph{bi-zero set} of \(\mu\).

\begin{theorem}[\cite{Jorgensen-Pedersen1998}]\label{J-P}
For a countable set \(\Lambda\subset\R\), define
\[
 Q_{\Lambda,\mu}(t)
 =\sum_{\lambda\in\Lambda}|\widehat\mu(t+\lambda)|^2,
 \qquad t\in\R.
\]
Then:
\begin{enumerate}[(i)]
\item \(\Lambda\) is a bi-zero set of \(\mu\) if and only if
\(Q_{\Lambda,\mu}(t)\leq1\) for every \(t\in\R\);
\item \(\Lambda\) is a spectrum of \(\mu\) if and only if
\(Q_{\Lambda,\mu}(t)=1\) for every \(t\in\R\).
\end{enumerate}
\end{theorem}

We also use the fact that \(Q_{\Lambda,\mu}\) is continuous whenever
\(\Lambda\) is a bi-zero set. Indeed, put \(e_t(x)=e^{2\pi itx}\). Bessel's
inequality for the orthonormal family \(E(-\Lambda)\) gives
\[
 \sum_{\lambda\in\Lambda}
 \bigl|\widehat\mu(t+\lambda)-\widehat\mu(s+\lambda)\bigr|^2
 \leq\|e_t-e_s\|_{L^2(\mu)}^2.
\]
Since \(Q_{\Lambda,\mu}\leq1\), the Cauchy--Schwarz inequality yields
\[
 |Q_{\Lambda,\mu}(t)-Q_{\Lambda,\mu}(s)|
 \leq2\|e_t-e_s\|_{L^2(\mu)}\longrightarrow0
 \qquad(s\to t).
\]

\begin{Lem}[\cite{Dai-He-Lau2014}]\label{Dai}
Let \(\mu=\nu*\omega\), and suppose that \(\Lambda\) is a bi-zero set of
\(\nu\). Then \(\Lambda\) is a bi-zero set of \(\mu\). If \(\omega\) is not
a Dirac measure, \(\Lambda\) cannot be a spectrum of \(\mu\).
\end{Lem}

\begin{Lem}\label{coefficient equiv}
Let \(c\in\R\setminus\{0\}\) and \(S_c(x)=cx\). The measure \(\mu\) is
spectral if and only if the pushforward \((S_c)_\#\mu\) is spectral.
\end{Lem}
\begin{proof}
Because \(\widehat{(S_c)_\#\mu}(t)=\widehat\mu(ct)\), a spectrum \(\Lambda\)
of \(\mu\) gives the spectrum \(c^{-1}\Lambda\) of \((S_c)_\#\mu\).
Applying the same argument to \(c^{-1}\) proves the converse.
\end{proof}

For \(n\in\mathbb Z\setminus\{0\}\), let \(v_2(n)\) be the largest
\(r\geq0\) such that \(2^r\mid n\), and put \(v_2(0)=\infty\). For
\(m/n\in\mathbb Q\), \(n\neq0\), define
\[
 v_2\!\left(\frac mn\right)=v_2(m)-v_2(n).
\]
We use the following characterization of Moran Bernoulli convolutions.

\begin{theorem}[\cite{Deng2023,Cao-Deng-Li-Wu2024}]\label{Bernoulli-sn}
Let \(b_k,d_k\in\mathbb Z\), with \(d_k\neq0\), and consider
\[
 \mu=\delta_{b_1^{-1}\{0,d_1\}}*
 \delta_{(b_1b_2)^{-1}\{0,d_2\}}*\cdots*
 \delta_{(b_1\cdots b_k)^{-1}\{0,d_k\}}*\cdots.
\]
Assume that \(|b_k|>|d_k|\) for \(k\geq2\) and that
\(\{|d_k|\}_{k\geq1}\) is bounded. Then \(\mu\) is spectral if and only if
the numbers
\[
 s_k=v_2\!\left(\frac{b_1b_2\cdots b_k}{2d_k}\right),
 \qquad k\geq1,
\]
are pairwise distinct.
\end{theorem}

\begin{Def}
Let \(R\in\mathbb Z\) satisfy \(|R|>1\), and let
\(B,L\subset\mathbb Z\) be finite sets with
\(\#B=\#L\geq2\). The triple \((R,B,L)\) is a \emph{Hadamard triple} if
\[
 \frac1{\sqrt{\#B}}
 \left(e^{-2\pi iR^{-1}b\ell}\right)_{b\in B,\,\ell\in L}
\]
is unitary.
\end{Def}

In particular, \((R,B,L)\) is a Hadamard triple precisely when \(L\) is a spectrum of
\(\delta_{R^{-1}B}\); see \cite{Dutkay-Han-Lai2019}.

\begin{Lem}[\cite{Li-Miao-Wang2024-1}]\label{finite Hadamard}
Let \(\{(N_j,B_j,L_j)\}_{j=1}^{n}\) be Hadamard triples. Set
\[
 \begin{aligned}
 N&=N_1N_2\cdots N_n,\\
 B&=(N_2\cdots N_n)B_1+\cdots+N_nB_{n-1}+B_n,\\
 L&=L_1+N_1L_2+\cdots+(N_1\cdots N_{n-1})L_n.
 \end{aligned}
\]
Then \((N,B,L)\) is a Hadamard triple.
\end{Lem}

Let \(\{(N_k,B_k,L_k)\}_{k\geq1}\) be a sequence of Hadamard triples.
Whenever the convolution exists, write
\begin{equation}\label{def-infinite convolution}
 \mu=\delta_{N_1^{-1}B_1}*\cdots*
 \delta_{(N_1\cdots N_k)^{-1}B_k}*\cdots.
\end{equation}
We use the decomposition \(\mu=\mu_k*\mu_{>k}\), where
\[
 \begin{aligned}
 \mu_k&=\delta_{N_1^{-1}B_1}*\cdots*
 \delta_{(N_1\cdots N_k)^{-1}B_k},\\
 \mu_{>k}&=\delta_{(N_1\cdots N_{k+1})^{-1}B_{k+1}}*
 \delta_{(N_1\cdots N_{k+2})^{-1}B_{k+2}}*\cdots.
 \end{aligned}
\]

\begin{Lem}[\cite{An-Fu-Lai2019}]\label{An}
Suppose that \(\{(N_k,B_k,L_k)\}_{k\geq1}\) is a sequence of Hadamard
triples satisfying
\[
 B_k\subset\{0,1,\ldots,N_k-1\}\quad(k\geq1),
 \qquad \liminf_{k\to\infty}\#B_k<\infty.
\]
Then the measure in \eqref{def-infinite convolution} is spectral and admits
a spectrum contained in \(\mathbb Z\).
\end{Lem}

The \emph{integral periodic zero set} of a Borel probability measure \(\mu\)
is
\[
 Z(\mu)=\{\xi\in\R:\widehat\mu(\xi+k)=0
 \text{ for every }k\in\mathbb Z\}.
\]
For \(E\subset\R\), let \(\mathcal P(E)\) denote the Borel probability
measures supported on \(E\). We equip \(\mathcal P(\R)\) with the usual weak
topology of probability measures: \(\eta_j\to\eta\) weakly if
\(\int f\,d\eta_j\to\int f\,d\eta\) for every \(f\in C_b(\R)\).
A family \(\Psi\subset\mathcal P(\R)\) is \emph{admissible} if
\(Z(\eta)=\varnothing\) for every \(\eta\in\overline{\Psi}\), where the
closure is taken in this topology.
It is \emph{tight} if, for every \(\varepsilon>0\), there is a compact set
\(K\subset\R\) such that \(\inf_{\eta\in\Psi}\eta(K)>1-\varepsilon\).

\begin{Lem}[\cite{An-Fu-Lai2019}]\label{empty}
If \(\mu\in\mathcal P([0,1])\), then \(Z(\mu)\neq\varnothing\) if and only if
\(\mu=\frac12(\delta_0+\delta_1)\).
\end{Lem}

For the decomposition following \eqref{def-infinite convolution}, define the
normalized tail
\[
 \omega_{>k}=(S_{N_1\cdots N_k})_\#\mu_{>k}
 =\delta_{N_{k+1}^{-1}B_{k+1}}*
 \delta_{(N_{k+1}N_{k+2})^{-1}B_{k+2}}*\cdots.
\]

\begin{theorem}[\cite{Li-Miao-Wang2024-1}]\label{Li}
Let \(\{(N_k,B_k,L_k)\}_{k\geq1}\) be Hadamard triples, and suppose that the
measure in \eqref{def-infinite convolution} exists. If some subsequence
\(\{\omega_{>n_j}\}_{j\geq1}\) is tight and admissible, then \(\mu\) is
spectral and has a spectrum contained in \(\mathbb Z\).
\end{theorem}

\begin{theorem}[\cite{Li-Miao-Wang2024-2}]\label{weak-conv}
Under the hypotheses of Theorem \ref{Li}, suppose that
\(\omega_{>n_j}\) converges weakly to \(\omega\) and
\(Z(\omega)=\varnothing\). Then \(\mu\) is spectral and has a spectrum
contained in \(\mathbb Z\).
\end{theorem}

We also need two consequences for sequences drawn from finitely many
Hadamard triples.

\begin{Lem}[\cite{Li-Miao-Wang2024-1}]\label{cl}
Let \(\{(N_j,B_j,L_j)\}_{j=1}^{m}\) be Hadamard triples and
\(\Sigma=\{1,\ldots,m\}^{\N}\). For
\(\sigma=(\sigma_1,\sigma_2,\ldots)\in\Sigma\), define
\[
 \mu_\sigma=\delta_{N_{\sigma_1}^{-1}B_{\sigma_1}}*\cdots*
 \delta_{(N_{\sigma_1}\cdots N_{\sigma_k})^{-1}B_{\sigma_k}}*\cdots.
\]
Then the family \(\Psi=\{\mu_\sigma:\sigma\in\Sigma\}\) is weakly closed.
\end{Lem}

\begin{Lem}[\cite{Li-Miao-Wang2024-2}]\label{finite-gcd}
Suppose that \(\{(N_k,B_k,L_k)\}_{k\geq1}\) is drawn from a finite family of
Hadamard triples. If
\[
 \gcd\left(\bigcup_{j=k}^{\infty}(B_j-B_j)\right)=1
 \qquad(k\geq1),
\]
then the measure in \eqref{def-infinite convolution} satisfies
\(Z(\mu)=\varnothing\).
\end{Lem}

Finally, we use the following elementary lemma.

\begin{Lem}[\cite{Deng2022}]\label{combined sum}
Let \(p_{i,j}>0\) and \(x_{i,j}\geq0\), where
\(1\leq i\leq m\) and \(1\leq j\leq n\). Assume that
\[
 \sum_{j=1}^{n}p_{i,j}=1\quad(1\leq i\leq m),
 \qquad
 \sum_{i=1}^{m}\max_{1\leq j\leq n}x_{i,j}\leq1.
\]
Then
\(\sum_{i=1}^{m}\sum_{j=1}^{n}p_{i,j}x_{i,j}=1
\)
if and only if
\(\sum_{i=1}^{m}x_{i,1}=1\) and  \(x_{i,1}=\cdots=x_{i,n}\) for \(1\leq i\leq m\).
\end{Lem}

\section{Existence and uniqueness of Moran measures}\label{sect.3}

In this section we construct the canonical pullback attractor and its
equal-weight probability measure, thereby proving Theorem \ref{thm-main1}.

\bigskip

\begin{proof}[Proof of Theorem \ref{thm-main1}]
For \(k\geq1\), recall that \(D_{n_r}=\{0,1,\ldots,n_r-1\}\) for
\(r\geq k\), and let
\[
 \Sigma_k=\prod_{r=k}^{\infty}D_{n_r}
\]
with the product topology. For
\(\sigma=(j_k,j_{k+1},\ldots)\in\Sigma_k\) and \(l\geq k\), set
\[
 \pi_{k,l}(\sigma)
 =\phi_{k,j_k}\circ\phi_{k+1,j_{k+1}}\circ\cdots
  \circ\phi_{l,j_l}(0).
\]
Each factor \(\phi_{r,j_r}\) has Lipschitz constant \(b_r^{-1}\), independently
of the sign \((-1)^{j_r}\). Repeated use of
\[
 |\phi_{r,j_r}(x)|\leq \frac{|x|+j_r}{b_r}
\]
therefore gives
\[
 \begin{aligned}
 |\pi_{k,l}(\sigma)|
 \leq\sum_{r=k}^{l}\frac{j_r}{b_kb_{k+1}\cdots b_r}\leq\sum_{r=k}^{l}\frac{n_r-1}{b_kb_{k+1}\cdots b_r}
 <\sum_{r=k}^{l}2^{-(r-k)}<2,
 \end{aligned}
\]
where we used \(n_r\leq b_r\) and \(b_s\geq2\). This estimate also applies
to a finite composition beginning at any later level.

If \(r>l\), the Lipschitz constant of
\(\phi_{k,j_k}\circ\cdots\circ\phi_{l,j_l}\) is
\((b_k\cdots b_l)^{-1}\). Put
\[
 y_{l+1,r}
 =\phi_{l+1,j_{l+1}}\circ\cdots\circ\phi_{r,j_r}(0).
\]
The preceding estimate gives \(|y_{l+1,r}|<2\), and hence
\[
 \begin{aligned}
 |\pi_{k,r}(\sigma)-\pi_{k,l}(\sigma)|
 &=\left|\phi_{k,j_k}\circ\cdots\circ\phi_{l,j_l}(y_{l+1,r})
 -\phi_{k,j_k}\circ\cdots\circ\phi_{l,j_l}(0)\right|\\
 &\leq\frac{|y_{l+1,r}|}{b_kb_{k+1}\cdots b_l}
 \leq\frac{2}{b_kb_{k+1}\cdots b_l}
 \leq 2^{k-l}.
 \end{aligned}
\]
The bound is independent of \(\sigma\), so \(\{\pi_{k,l}\}_{l\geq k}\)
is uniformly Cauchy. Each \(\pi_{k,l}\) depends on only finitely many
coordinates and is therefore continuous in the product topology. Consequently,
it converges uniformly to a continuous map
\(\pi_k:\Sigma_k\to[-2,2]\).

Each factor \(D_{n_r}\) is finite and hence compact. Thus \(\Sigma_k\) is
compact by Tychonoff's theorem. Define \(E_k=\pi_k(\Sigma_k)\).
It follows that \(E_k\) is a nonempty compact subset of \([-2,2]\), so the
family \(\{E_k\}_{k\geq1}\) is uniformly bounded.

For \(i\in D_{n_k}\) and
\(\sigma=(j_{k+1},j_{k+2},\ldots)\in\Sigma_{k+1}\), let
\(i\sigma=(i,j_{k+1},j_{k+2},\ldots)\). At every finite level,
\[
 \pi_{k,l}(i\sigma)=\phi_{k,i}\bigl(\pi_{k+1,l}(\sigma)\bigr).
\]
Passing to the limit and using the continuity of \(\phi_{k,i}\) gives
\[
 \pi_k(i\sigma)=\phi_{k,i}\bigl(\pi_{k+1}(\sigma)\bigr),
 \qquad i\in D_{n_k},\quad \sigma\in\Sigma_{k+1},
\]
and hence
\[
 E_k=\bigcup_{i\in D_{n_k}}\phi_{k,i}(E_{k+1})
 =\Phi_k(E_{k+1}).
\]

For \(a_{l+1}\in E_{l+1}\subset[-2,2]\), the same Lipschitz estimate yields
\[
 \left|
 \phi_{k,j_k}\circ\cdots\circ\phi_{l,j_l}(a_{l+1})
 -\phi_{k,j_k}\circ\cdots\circ\phi_{l,j_l}(0)
 \right|
 \leq\frac{2}{b_k\cdots b_l}\leq2^{k-l}.
\]
The right-hand side is independent of the digits and of \(a_{l+1}\), and
tends to zero as \(l\to\infty\). This proves precisely the uniform convergence
asserted in \eqref{conv-point}.

To prove uniqueness of the compact family, suppose that
\(\{F_k\}_{k\geq1}\) is another uniformly bounded family of nonempty compact
sets satisfying \eqref{equi-supp}, and let
\[
 M=\sup_{k\geq1}\sup_{x\in F_k}|x|<\infty.
\]
Fix \(x\in F_k\). Since \(F_k=\Phi_k(F_{k+1})\), there exist
\(j_k\in D_{n_k}\) and \(x_{k+1}\in F_{k+1}\) such that
\(x=\phi_{k,j_k}(x_{k+1})\). Repeating this choice gives a single sequence
\(\sigma=(j_k,j_{k+1},\ldots)\in\Sigma_k\) and points
\(x_{l+1}\in F_{l+1}\) satisfying
\[
 x=\phi_{k,j_k}\circ\cdots\circ\phi_{l,j_l}(x_{l+1})
 \qquad(l\geq k).
\]
Therefore
\[
 |x-\pi_{k,l}(\sigma)|
 \leq \frac{|x_{l+1}|}{b_k\cdots b_l}
 \leq \frac{M}{b_k\cdots b_l}\longrightarrow0.
\]
Since \(\pi_{k,l}(\sigma)\to\pi_k(\sigma)\), we obtain
\(x=\pi_k(\sigma)\in E_k\). Thus \(F_k\subset E_k\).

Conversely, fix \(\sigma=(j_k,j_{k+1},\ldots)\in\Sigma_k\). For every
\(l\geq k\), choose any \(y_{l+1}\in F_{l+1}\). Iterating
\eqref{equi-supp} shows that
\[
 z_l:=\phi_{k,j_k}\circ\cdots\circ\phi_{l,j_l}(y_{l+1})\in F_k.
\]
Moreover,
\[
 |z_l-\pi_{k,l}(\sigma)|
 \leq \frac{M}{b_k\cdots b_l}\longrightarrow0,
\]
and \(\pi_{k,l}(\sigma)\to\pi_k(\sigma)\). Hence
\(z_l\to\pi_k(\sigma)\). Since \(F_k\) is closed,
\(\pi_k(\sigma)\in F_k\), proving \(E_k\subset F_k\). Therefore
\(E_k=F_k\) for every \(k\), and part (i) follows.

For (ii), let
\[
 \omega_k=\bigotimes_{r=k}^{\infty}
 \left(\frac1{n_r}\sum_{j=0}^{n_r-1}\delta_j\right)
\]
be the product probability measure on \(\Sigma_k\), and set
\(\mu_k=(\pi_k)_\#\omega_k\). The product measure has full support, so
\(\supp\omega_k=\Sigma_k\): indeed, every nonempty cylinder has positive
\(\omega_k\)-measure. Since \(\pi_k\) is continuous and \(\Sigma_k\) is
compact,
\[
 \supp\mu_k=\pi_k(\supp\omega_k)=\pi_k(\Sigma_k)=E_k.
\]
To justify the middle equality, let \(U\) be a neighborhood of
\(\pi_k(\sigma)\). By continuity, \(\pi_k^{-1}(U)\) contains a neighborhood
of \(\sigma\), which has positive \(\omega_k\)-measure. Thus every point of
\(E_k\) belongs to \(\supp\mu_k\); the reverse inclusion is immediate.

We next verify the measure recursion. If \(f\in C_b(\R)\), then conditioning
the product measure on its first coordinate and using the coding identity gives
\[
 \begin{aligned}
 \int f\,d\mu_k
 &=\int_{\Sigma_k}f(\pi_k(\sigma))\,d\omega_k(\sigma)\\
 &=\frac1{n_k}\sum_{i=0}^{n_k-1}
   \int_{\Sigma_{k+1}}
   f\bigl(\phi_{k,i}(\pi_{k+1}(\sigma))\bigr)
   \,d\omega_{k+1}(\sigma)\\
 &=\frac1{n_k}\sum_{i=0}^{n_k-1}
   \int f\,d\bigl((\phi_{k,i})_\#\mu_{k+1}\bigr).
 \end{aligned}
\]
Since this holds for every bounded continuous \(f\),
\[
 \mu_k=\frac1{n_k}\sum_{i=0}^{n_k-1}
 (\phi_{k,i})_\#\mu_{k+1}.
\]

It remains to prove uniqueness. On the probability measures supported in
\([-2,2]\), the Kantorovich--Rubinstein dual formula for the
\(1\)-Wasserstein distance is
\[
 W_1(\eta,\eta')
 =\sup_{\operatorname{Lip}(f)\leq1}
 \left|\int f\,d\eta-\int f\,d\eta'\right|.
\]
Define the level-\(k\) averaging operator by
\[
 T_k\eta=\frac1{n_k}\sum_{i=0}^{n_k-1}(\phi_{k,i})_\#\eta,
\qquad \eta\in\mathcal P([-2,2]).
\]
If \(f\) is \(1\)-Lipschitz, then
\(f\circ\phi_{k,i}\) is \(b_k^{-1}\)-Lipschitz. Therefore
\[
 \begin{aligned}
 &\left|\int f\,d(T_k\eta)-\int f\,d(T_k\eta')\right|\\
 &\quad=\left|\frac1{n_k}\sum_{i=0}^{n_k-1}
 \left(\int f\circ\phi_{k,i}\,d\eta
       -\int f\circ\phi_{k,i}\,d\eta'\right)\right|\\
 &\quad\leq\frac1{n_k}\sum_{i=0}^{n_k-1}
 \frac1{b_k}W_1(\eta,\eta')
 =\frac1{b_k}W_1(\eta,\eta').
 \end{aligned}
\]
Taking the supremum over all such \(f\) proves the contraction estimate
\[
 W_1(T_k\eta,T_k\eta')\leq b_k^{-1}W_1(\eta,\eta').
\]

Consequently, if \(\{\mu_k'\}_{k\geq1}\) is another family satisfying
\eqref{equi-mea} and supported on \(E_k\), then for every \(m\geq1\),
\[
 \begin{aligned}
 W_1(\mu_k,\mu_k')
 \leq\frac{W_1(\mu_{k+m},\mu_{k+m}')} {b_kb_{k+1}\cdots b_{k+m-1}}\leq\frac{4}{b_kb_{k+1}\cdots b_{k+m-1}} \leq4\cdot2^{-m}.
 \end{aligned}
\]
Here the second inequality follows from
\(W_1(\eta,\eta')\leq4\) for probability measures supported in
\([-2,2]\). The last bound tends to zero as
\(m\to\infty\), so \(W_1(\mu_k,\mu_k')=0\). Hence
\(\mu_k=\mu_k'\) for every \(k\), completing the proof.
\end{proof}

\section{Proof of Theorem \ref{thm-main2}} \label{sect.4}
In this section we study the spectrality of $\mu:=\mu_1$. We assume throughout
that $n_k=2p_k$, where $p_k\in\N$, and write
$B_0=1$ and $B_k=b_1\cdots b_k$ for $k\geq1$. Then
\begin{equation}\label{def-mu_k}
 \mu_k=\frac1{2p_k}\sum_{j=0}^{2p_k-1}
 (\phi_{k,j})_\#\mu_{k+1},\qquad k\geq1.
\end{equation}

\begin{Lem}\label{lem-Fourier transform}
For \(k\geq1\), define
\[
 f_k(u)=\frac1{p_k}\sum_{j=0}^{p_k-1}
 \cos\bigl((4j+1-b_{k+1}^{-1})\pi u\bigr).
\]
Then
\[
 \widehat\mu(t)=e^{-\pi it/b_1}
 \prod_{k=1}^{\infty}f_k(t/B_k),\qquad t\in\R,
\]
and the product converges locally uniformly on \(\R\).
\end{Lem}

\begin{proof}
Taking Fourier transforms in \eqref{def-mu_k} gives
\[
 \widehat\mu_k(t)
 =\frac1{2p_k}\sum_{j=0}^{2p_k-1}
 e^{2\pi it(-1)^j j/b_k}
 \widehat\mu_{k+1}\bigl((-1)^jt/b_k\bigr).
\]
Put
\[
 v_k(t)=
 \begin{pmatrix}
  \operatorname{Re}\widehat\mu_k(t)\\
  \operatorname{Im}\widehat\mu_k(t)
 \end{pmatrix}.
\]
Separating the real and imaginary parts yields
\begin{equation}\label{eq-iteration of transform}
 v_k(t)=M_k(t/b_k)v_{k+1}(t/b_k).
\end{equation}
Write
\[
 M_k(u)=
 \begin{pmatrix}
  a_{k,1}(u)&a_{k,2}(u)\\
  a_{k,3}(u)&a_{k,4}(u)
 \end{pmatrix},
\]
where
\[
 \begin{aligned}
 a_{k,1}(u)&=\frac1{2p_k}\sum_{j=0}^{2p_k-1}\cos(2\pi ju),&
 a_{k,2}(u)&=-\frac1{2p_k}\sum_{j=0}^{2p_k-1}\sin(2\pi ju),\\
 a_{k,3}(u)&=\frac1{2p_k}\sum_{j=0}^{2p_k-1}(-1)^j\sin(2\pi ju),&
 a_{k,4}(u)&=\frac1{2p_k}\sum_{j=0}^{2p_k-1}(-1)^j\cos(2\pi ju).
 \end{aligned}
\]
Pairing the terms with indices \(2j\) and \(2j+1\) gives the rank-one
factorization
\[
 M_k(u)=\alpha(u)^{\mathsf T}\beta_k(u),
\]
where
\[
 \begin{aligned}
 \alpha(u)&=(\cos\pi u,-\sin\pi u),\\
 \beta_k(u)&=\frac1{p_k}\left(
  \sum_{j=0}^{p_k-1}\cos((4j+1)\pi u),
  -\sum_{j=0}^{p_k-1}\sin((4j+1)\pi u)\right).
 \end{aligned}
\]
The cosine addition formula therefore gives
\begin{equation}\label{rank-one-scalar}
 \beta_k(u)\alpha(u/b_{k+1})^{\mathsf T}
 =f_k(u).
\end{equation}

Iterating \eqref{eq-iteration of transform} and using
\eqref{rank-one-scalar}, we obtain, for \(N\geq2\),
\begin{equation}\label{iterated-matrix-product}
 v_1(t)=\alpha(t/B_1)^{\mathsf T}
 \left(\prod_{k=1}^{N-1}f_k(t/B_k)\right)
 \beta_N(t/B_N)v_{N+1}(t/B_N).
\end{equation}
The measures \(\mu_{N+1}\) are uniformly supported in \([-2,2]\), so
\(v_{N+1}(t/B_N)\to(1,0)^{\mathsf T}\), locally uniformly in \(t\).
Furthermore,
\[
\frac{p_N}{B_N}\leq\frac{1}{2B_{N-1}}\longrightarrow0.
\]
Hence, for every \(M>0\),
\[
\max_{\substack{|t|\leq M\\0\leq j<p_N}}
\left|\frac{(4j+1)\pi t}{B_N}\right|
\leq 4\pi M\frac{p_N}{B_N}\longrightarrow0.
\]
It follows directly from the definition of \(\beta_N\) that \(\beta_N(t/B_N)\to(1,0)\), also locally uniformly.
Thus the last scalar factor in \eqref{iterated-matrix-product} tends to \(1\).

The infinite product is locally uniformly convergent. Indeed,
\(p_k/B_k\leq2^{-k}\), and for \(|t|\leq M\),
\[
 \begin{aligned}
 0\leq1-f_k(t/B_k)
 &\leq\frac{\pi^2t^2}{2p_kB_k^2}
 \sum_{j=0}^{p_k-1}(4j+1-b_{k+1}^{-1})^2\\
 &\leq8\pi^2M^2\frac{p_k^2}{B_k^2}
 \leq8\pi^2M^2\,4^{-k}.
 \end{aligned}
\]
Letting \(N\to\infty\) in \eqref{iterated-matrix-product} now gives
\[
 v_1(t)=\alpha(t/B_1)^{\mathsf T}
 \prod_{k=1}^{\infty}f_k(t/B_k).
\]
Since \(\alpha(t/B_1)^{\mathsf T}\) consists of the real and imaginary parts
of \(e^{-\pi it/b_1}\), the asserted formula follows.
\end{proof}

For each \(k\geq1\), set
\begin{equation}\label{eq-D_k}
 \mathcal D_k
 =D_{p_k}\oplus
 \left(p_k-\frac{1+b_{k+1}^{-1}}2\right)D_2.
\end{equation}
Consider the infinite convolution
\begin{equation}\label{eq-infinite convolution-equivalent}
 \nu=\delta_{B_1^{-1}\mathcal D_1}*
 \delta_{B_2^{-1}\mathcal D_2}*\cdots*
 \delta_{B_k^{-1}\mathcal D_k}*\cdots.
\end{equation}

\begin{Prop}\label{equivalent}
The Moran measure \(\mu=\mu_1\) is spectral if and only if the measure
\(\nu\) in \eqref{eq-infinite convolution-equivalent} is spectral.
\end{Prop}

\begin{proof}
 Put
\(c_k=p_k-(1+b_{k+1}^{-1})/2\).
Since \(\mathcal D_k=D_{p_k}\oplus c_kD_2\), its mask polynomial factors as
\(m_{\mathcal D_k}=m_{D_{p_k}}m_{c_kD_2}\). For \(u\notin\mathbb Z\),
\[
 \begin{aligned}
 m_{D_{p_k}}(u)
 &=e^{(p_k-1)\pi iu}
   \frac{\sin(p_k\pi u)}{p_k\sin(\pi u)},\\
 m_{c_kD_2}(u)
 &=\frac{1+e^{2\pi ic_ku}}2
   =e^{\pi ic_ku}\cos(\pi c_ku).
 \end{aligned}
\]
The finite cosine-sum formula also gives
\[
 f_k(u/2)=\frac{\sin(p_k\pi u)}{p_k\sin(\pi u)}\cos(\pi c_ku).
\]
Multiplying the two mask identities gives
\begin{equation}\label{eq-mask polyn.}
	m_{\mathcal D_k}(u)
	=f_k(u/2)\exp\!\left[
	\pi i\left(2p_k-\frac{3+b_{k+1}^{-1}}2\right)u\right].
\end{equation}
Continuity then extends it to every \(u\in\R\).

Applying \eqref{eq-mask polyn.} at \(u=t/B_k\) and using Lemma
\ref{lem-Fourier transform}, we obtain
\[
 \widehat\nu(t)=e^{\pi i\theta t}\widehat\mu(t/2),
\]
where
\[
 \theta=\frac1{2b_1}
 +\sum_{k=1}^{\infty}
 \left(2p_k-\frac{3+b_{k+1}^{-1}}2\right)\frac1{B_k}.
\]
The series converges because \(2p_k\leq b_k\) and \(B_{k-1}\geq2^{k-1}\).
The phase has modulus one; hence, for every countable
\(\Lambda\subset\R\),
\[
 Q_{\Lambda,\nu}(t)
 =\sum_{\lambda\in\Lambda}
 \left|\widehat\mu\!\left(\frac{t+\lambda}{2}\right)\right|^2
 =Q_{\Lambda/2,\mu}(t/2).
\]
The conclusion follows from Theorem \ref{J-P}.
\end{proof}

By Proposition \ref{equivalent}, it remains to study \(\nu\). The Bernoulli
case \(p_k=1\) for all \(k\) also shows that the sufficient conditions in
Theorem \ref{thm-main2} need not be necessary without the assumption
\(4\mid n_k\).

\begin{theorem}\label{Bernoulli situation}
Assume that \(p_k=1\) for every \(k\geq1\) and that
\(\{b_k\}_{k\geq1}\) is bounded. If \(2\mid b_k\) for every \(k\geq2\),
then the measure \(\nu\) in \eqref{eq-infinite convolution-equivalent} is
spectral. The converse fails: there are spectral examples for which some
\(b_k\) is odd.
\end{theorem}

\begin{proof}
When \(p_k=1\),
\[
 \mathcal D_k=\left\{0,\frac{1-b_{k+1}^{-1}}2\right\}.
\]
With \(\widetilde D_k=2b_{k+1}\mathcal D_k
=\{0,b_{k+1}-1\}\), we can rewrite \(\nu\) as
\[
 \nu=\delta_{(2b_1b_2)^{-1}\widetilde D_1}*
 \delta_{(2b_1b_2b_3)^{-1}\widetilde D_2}*
 \delta_{(2b_1b_2b_3b_4)^{-1}\widetilde D_3}*\cdots.
\]
This is a Moran Bernoulli convolution to which Theorem
\ref{Bernoulli-sn} applies. If every \(b_k\), \(k\geq2\), is even, then
\[
 \begin{aligned}
 s_{k+1} =v_2\!\left(\frac{b_1\cdots b_{k+2}}{b_{k+2}-1}\right)
  =v_2(b_1\cdots b_{k+2})>v_2(b_1\cdots b_{k+1})
  =v_2\!\left(\frac{b_1\cdots b_{k+1}}{b_{k+1}-1}\right)
  =s_k .
 \end{aligned}
\]
Thus the \(s_k\) are pairwise distinct, and \(\nu\) is spectral.

For the converse assertion, take \(b_3=7\) and \(b_k=8\) for \(k\neq3\).
Then
\[
 s_1=6,\qquad s_2=5,\qquad s_k=3k\quad(k\geq3).
\]
These numbers are pairwise distinct, so Theorem \ref{Bernoulli-sn} again
implies that \(\nu\) is spectral, although \(b_3\) is odd.
\end{proof}

For general \(p_k\), put
\[
 a_k=b_{k+1}(2p_k-1)-1,
 \qquad
 \eta_k=\delta_{B_k^{-1}D_{p_k}},
 \qquad
 \tau_k=\delta_{(2B_{k+1})^{-1}a_kD_2}.
\]
Thus \(\eta_k\) is the consecutive-digit factor, while \(\tau_k\) is the
two-point factor.
Equation \eqref{eq-D_k} gives
\[
 \delta_{B_k^{-1}\mathcal D_k}=\eta_k*\tau_k,
 \qquad \nu=(\eta_1*\tau_1)*(\eta_2*\tau_2)*\cdots.
\]
Since convolution is commutative, we may use the ordering
\begin{equation}\label{eq-form of measure nu}
 \nu=\eta_1*\eta_2*\tau_1*\eta_3*\tau_2*\eta_4*\tau_3*\cdots.
\end{equation}

We remove the trivial factors $\eta_k=\delta_0$
for which $p_k=1$ and relabel the remaining factors. The first terms are
defined as follows:
\begin{enumerate}[(i)]
\item if $p_1=p_2=1$, let \( b'_1=2B_2, \ D'_1=a_1D_2\);
\item if $p_1=1$ and $p_2\neq1$, let
\(b'_1=B_2, \ D'_1=D_{p_2}, \  b'_2=2, \ D'_2=a_1D_2\);
\item if $p_1\neq1$ and $p_2=1$, let
\(b'_1=b_1, \ D'_1=D_{p_1}, \   b'_2=2b_2, \ D'_2=a_1D_2\);
\item if $p_1\neq1$ and $p_2\neq1$, let \(b'_1=b_1, \  D'_1=D_{p_1}, \  b'_2=b_2, \  D'_2=D_{p_2}, \  b'_3=2, \ D'_3=a_1D_2\).
\end{enumerate}

For $k\geq3$, put
\[
 m_k=2k-2-\#\{n<k:p_n=1\}.
\]
The remaining terms are relabeled by the following two rules:
\medskip

\noindent\textup{(v)} If $p_k\neq1$, let  \(b'_{m_k}=\frac{b_k}{2}, \ D'_{m_k}=D_{p_k}, \  b'_{m_k+1}=2, \ D'_{m_k+1}=a_{k-1}D_2\);

\noindent\textup{(vi)} If $p_k=1$, let \(b'_{m_k}=b_k, \ D'_{m_k}=a_{k-1}D_2\).

The initial rules (i)--(iv) ensure that the product of the scale factors
introduced before the third block is \(2B_2\). We now verify the remaining
rules inductively. Suppose that, before the \(k\)th block, the product of the
previously introduced scale factors is \(2B_{k-1}\).

If \(p_k\neq1\), rule (v) introduces the two scale factors \(b_k/2\) and
\(2\). Their successive cumulative products are
\[
(2B_{k-1})\frac{b_k}{2}=B_k,
\qquad
(2B_{k-1})\frac{b_k}{2}\cdot2=2B_k.
\]
Together with the corresponding digit sets \(D_{p_k}\) and
\(a_{k-1}D_2\), these give precisely
\[
\delta_{B_k^{-1}D_{p_k}}=\eta_k
\quad\text{and}\quad
\delta_{(2B_k)^{-1}a_{k-1}D_2}=\tau_{k-1}.
\]

If \(p_k=1\), then \(\eta_k=\delta_0\) is omitted. Rule (vi) introduces only
the scale factor \(b_k\), and
\[
(2B_{k-1})b_k=2B_k,
\]
which gives the remaining factor
\[
\delta_{(2B_k)^{-1}a_{k-1}D_2}=\tau_{k-1}.
\]
In either case, the cumulative product after the \(k\)th block is \(2B_k\),
so the induction continues. Hence the relabeled convolution contains exactly
the nontrivial factors in \eqref{eq-form of measure nu}, in the same order.
That is,
\begin{equation}\label{rearrangement}
\nu=\delta_{{b'_1}^{-1}D'_1}*\delta_{(b'_1 b'_2)^{-1}D'_2}*\cdots
\end{equation}
with $\#D'_k\geq2$ for every $k$.

For later use, decompose \(\nu\) as
\[
 \nu=\nu_k*\nu_{>k}
\]
where
\[
 \begin{aligned}
 \nu_k&=\delta_{{b'_1}^{-1}D'_1}*
 \delta_{(b'_1b'_2)^{-1}D'_2}*\cdots*
 \delta_{(b'_1\cdots b'_k)^{-1}D'_k},\\
 \nu_{>k}&=\delta_{(b'_1\cdots b'_{k+1})^{-1}D'_{k+1}}*
 \delta_{(b'_1\cdots b'_{k+2})^{-1}D'_{k+2}}*\cdots.
 \end{aligned}
\]
We normalize the tail by the dilation introduced in Lemma
\ref{coefficient equiv}:
\begin{equation}\label{omega.n}
\omega_{>k}:=(S_{b'_1\cdots b'_k})_\#\nu_{>k}
=\delta_{{b'_{k+1}}^{-1}D'_{k+1}}*
\delta_{(b'_{k+1}b'_{k+2})^{-1}D'_{k+2}}*\cdots.
\end{equation}

\begin{Prop}\label{4.4}
Assume that \(p_1\mid b_1, \  p_2\mid b_2, \  2\mid b_2\), and  \(2p_k\mid b_k\) for \(k\geq3\). Then there exists a sequence of finite sets
\(\{L_k\}_{k\geq1}\), with \(L_k\subset\mathbb Z\), such that
\((b'_k,D'_k,L_k)\) is a Hadamard triple for every \(k\geq1\).
\end{Prop}
\begin{proof}
We first record the two elementary Hadamard triples used below. If \(m\mid R\),
then \(\left(R,D_m,\frac{R}{m}D_m\right)\) is a Hadamard triple, since its defining matrix is
\[
 \frac1{\sqrt m}
 \left(e^{-2\pi i rs/m}\right)_{0\leq r,s<m}.
\]
If \(R\) is even and \(a\) is odd, then \(\left(R,aD_2,\frac{R}{2}D_2\right)\) is also a Hadamard triple. It remains to verify that every pair \((b'_j,D'_j)\) defined by rules
(i)--(vi) is covered by one of the two cases above.

In the initial rules (i)--(iv), a nontrivial consecutive-digit factor \(D_{p_1}\) has scale
\(b_1\), while a nontrivial factor \(D_{p_2}\) has scale \(b_2\) or \(B_2\).
The required divisibility follows from \(p_1\mid b_1,  p_2\mid b_2\). The initial two-point factor has digit set \(a_1D_2\) and an even scale. Moreover, \(a_1=b_2(2p_1-1)-1\) is odd because \(b_2\) is even.

For \(k\geq3\), suppose first that \(p_k>1\). Rule (v) produces a
consecutive-digit factor \(D_{p_k}\) with scale \(b_k/2\) and a two-point
factor \(a_{k-1}D_2\) with scale \(2\). Since \(2p_k\mid b_k\), we have \(p_k\mid\frac{b_k}{2}\) and
\(a_{k-1}=b_k(2p_{k-1}-1)-1\) is odd. If \(p_k=1\), the consecutive-digit factor is trivial and is omitted.
Rule (vi) then leaves only the two-point factor \(a_{k-1}D_2\) with scale
\(b_k\). In this case \(b_k\) is even and \(a_{k-1}\) is again odd.

For each consecutive-digit factor \(D'_j=D_m\), define \(L_j=\frac{b'_j}{m}D_m\), and for each two-point factor define
\(L_j=\frac{b'_j}{2}D_2\).  The divisibility and parity properties verified above ensure that
\(L_j\subset\mathbb Z\), and \((b'_j,D'_j,L_j)\) is a Hadamard triple for every \(j\geq1\).
\end{proof}

\begin{theorem}\label{thm3- spectral measure}
Assume that \(\{b_k\}_{k\geq1}\) is bounded and that \(p_2\mid b_2, \ 2\mid b_2, \ 2p_k\mid b_k\) for \(k\geq3\). Then the measure \(\nu\) in \eqref{rearrangement} is spectral.
\end{theorem}
\begin{proof}
By Lemma \ref{coefficient equiv}, the value of \(b_1\) affects \(\nu\) only
through a global dilation. We may therefore assume that \(p_1\mid b_1\).

We distinguish three cases according to the number of indices for which
\(p_k\geq2\).

\emph{Case 1.} If \(p_k=1\) for every \(k\geq1\), the result follows from
Theorem \ref{Bernoulli situation}.

\emph{Case 2.} Suppose that \(p_k\geq2\) for only finitely many \(k\).
Then there exists \(\kappa\ge 3\) such that \(p_k=1\) for \(k\geq\kappa\). Then \(m_{\kappa+j}=m_\kappa+j\) for every \(j\geq1\). Let \(\{L_k\}_{k\geq1}\) be given by Proposition \ref{4.4}, and set
\[
 \begin{aligned}
 B&=b'_1b'_2\cdots b'_{m_\kappa},\\
 L&=L_1+b'_1L_2+\cdots+
 (b'_1b'_2\cdots b'_{m_\kappa-1})L_{m_\kappa},\\
 D&=(b'_2\cdots b'_{m_\kappa})D'_1+
 (b'_3\cdots b'_{m_\kappa})D'_2+\cdots+
 b'_{m_\kappa}D'_{m_\kappa-1}+D'_{m_\kappa}.
 \end{aligned}
\]
Then
\[
 \nu_{m_\kappa}
 =\delta_{{b'_1}^{-1}D'_1}*
 \delta_{(b'_1b'_2)^{-1}D'_2}*\cdots*
 \delta_{(b'_1\cdots b'_{m_\kappa})^{-1}D'_{m_\kappa}}
 =\delta_{B^{-1}D}.
\]
By Lemma \ref{finite Hadamard},
\((B,D,L)\) is a Hadamard triple. Hence \(L\) is a spectrum of
\(\nu_{m_\kappa}\), and
\(\supp\nu_{m_\kappa}\subset B^{-1}\mathbb Z\).

Since \(p_k=1\) for \(k\geq\kappa\), the normalized tail in
\eqref{omega.n} has the form
\[
 \begin{aligned}
 \omega_{>m_\kappa}
 &=\delta_{b_{\kappa+1}^{-1}(b_{\kappa+1}-1)D_2}*
 \delta_{(b_{\kappa+1}b_{\kappa+2})^{-1}
 (b_{\kappa+2}-1)D_2}*\cdots .
 \end{aligned}
\]
Thus its defining triples are
\[
\left(
b_{\kappa+j},
(b_{\kappa+j}-1)D_2,
\frac{b_{\kappa+j}}2D_2
\right),
\qquad j\geq1.
\]
They satisfy the hypotheses of Lemma \ref{An}; hence
\(\omega_{>m_\kappa}\) has a spectrum
\(\Gamma\subset\mathbb Z\). It follows from the definition of the normalized
tail that \(B\Gamma\) is a spectrum of \(\nu_{>m_\kappa}\).

Set
\[
 \Lambda=L\oplus B\Gamma.
\]
The sum is disjoint because the elements of \(L\) are distinct modulo \(B\).
Also,
\(\widehat{\nu_{m_\kappa}}(t+B\gamma)
=\widehat{\nu_{m_\kappa}}(t)\) for \(\gamma\in\mathbb Z\). Therefore,
for every \(t\in\R\),
\[
 \begin{aligned}
 Q_{\Lambda,\nu}(t)
 &=\sum_{\lambda\in L}\sum_{\gamma\in\Gamma}
 |\widehat{\nu_{m_\kappa}}(t+\lambda+B\gamma)|^2
 |\widehat{\nu_{>m_\kappa}}(t+\lambda+B\gamma)|^2\\
 &=\sum_{\lambda\in L}|\widehat{\nu_{m_\kappa}}(t+\lambda)|^2
 \sum_{\gamma\in\Gamma}
 |\widehat{\nu_{>m_\kappa}}(t+\lambda+B\gamma)|^2\\
 &=\sum_{\lambda\in L}
 |\widehat{\nu_{m_\kappa}}(t+\lambda)|^2=1.
 \end{aligned}
\]
Theorem \ref{J-P} shows that \(\Lambda\) is a spectrum of \(\nu\).

\emph{Case 3.} Suppose that \(p_k\geq2\) for infinitely many \(k\).
For every such \(k\), rule (v) gives \(D'_{m_k}=D_{p_k}\).  Since \(p_k\geq2\), the difference set \(D_{p_k}-D_{p_k}\) contains \(1\),
and hence
\[
\gcd(D'_{m_k}-D'_{m_k})=1.
\]
Moreover, \(m_k =2k-2-\#\{n<k:p_n=1\} \geq k-1\). 
Thus the indices \(m_k\) for which \(p_k\geq2\) are unbounded. Consequently,
for every \(r\geq1\),
\begin{equation}\label{tail-gcd}
	\gcd\left(
	\bigcup_{j=r}^{\infty}(D'_j-D'_j)
	\right)=1.
\end{equation}

Since \(\{b_k\}_{k\geq1}\) is bounded, the triples
\(\{(b'_k,D'_k,L_k)\}_{k\geq1}\) are drawn from a finite family of
Hadamard triples. For each \(l\geq1\), the defining sequence of
\(\omega_{>l}\) is
\[\{(b'_{l+j},D'_{l+j},L_{l+j})\}_{j\geq1}.\]
By \eqref{tail-gcd},
\[
\gcd\left(
\bigcup_{j=k}^{\infty}
(D'_{l+j}-D'_{l+j})
\right)
=
\gcd\left(
\bigcup_{r=l+k}^{\infty}
(D'_r-D'_r)
\right)
=1
\]
for every \(k\geq1\). Lemma \ref{finite-gcd} therefore yields
\(Z(\omega_{>l})=\varnothing\) for \(l\geq1 \).

We now consider weak limits of the normalized tails. Since the defining
triples range over a finite family, a diagonal argument gives a sequence
\(l_j\to\infty\) such that the defining sequences of
\(\omega_{>l_j}\) converge coordinatewise. The digit sets are uniformly
bounded and every scale factor is at least \(2\). Consequently,   all these
tails are supported on a common compact set, and their contributions beyond the first \(N\) factors are bounded
uniformly by \(C2^{-N}\) for some constant $C>0$. It follows that \(\omega_{>l_j}\to\omega'\)
weakly for an infinite convolution \(\omega'\) generated by the same finite
family of Hadamard triples. Moreover, for every fixed \(M\geq0\),
\(\omega_{>l_j+M}\to\omega'_{>M}\) weakly.

Let \(\Psi\) denote the family of infinite convolutions generated by this
finite collection of Hadamard triples, and let \(\Psi_0\subset\Psi\) consist
of those measures whose defining sequences contain infinitely many
consecutive digit sets \(D_p\) with \(p\geq2\). By Lemma \ref{cl},
\(\omega'\in\Psi\).

\textup{(i)}
Suppose that \(\omega'\in\Psi_0\). Then every tail of its defining sequence
contains a digit set \(D_p\) with \(p\geq2\). Since \(D_p-D_p\) contains
\(1\), the gcd condition in Lemma \ref{finite-gcd} holds for every tail.
Therefore \(Z(\omega')=\varnothing\).  Theorem \ref{weak-conv}, applied to
\(\omega_{>l_j}\to\omega'\), shows that \(\nu\) is spectral.

\textup{(ii)}
Suppose that \(\omega'\in\Psi\setminus\Psi_0\). Then its defining sequence
contains only finitely many consecutive digit sets \(D_p\) with \(p\geq2\).
By rules (v)--(vi), after deleting sufficiently many initial factors, all
the remaining factors are of Bernoulli type. Thus there exists \(M\geq0\)
such that, writing the successive scale factors of \(\omega'_{>M}\) as
\(c_1,c_2,\ldots\), we have
\[
\omega'_{>M}
=
\delta_{c_1^{-1}(c_1-1)D_2}
*
\delta_{(c_1c_2)^{-1}(c_2-1)D_2}
*\cdots .
\]
The largest point in the support of its \(n\)-factor approximation is
\[
\sum_{r=1}^{n}
\frac{c_r-1}{c_1\cdots c_r}
=
1-\frac1{c_1\cdots c_n}
\longrightarrow1.
\]
Hence \(\supp\omega'_{>M}\subset[0,1]\). The point \(1\) can be obtained only by choosing the maximal digit at every
level. Since each of the two digits has probability \(1/2\),
\[
\omega'_{>M}(\{1\})
\leq 2^{-n}
\qquad(n\geq1),
\]
and therefore \(\omega'_{>M}(\{1\})=0\). In particular, \(\omega'_{>M}\neq\frac12(\delta_0+\delta_1)\). Lemma \ref{empty} now gives
\(Z(\omega'_{>M})=\varnothing\). Finally, Theorem \ref{weak-conv}, applied to \(\omega_{>l_j+M}\to \omega'_{>M}\),
shows that \(\nu\) is spectral.
\end{proof}

\bigskip

\begin{proof}[Proof of Theorem \ref{thm-main2}]
The result follows from Proposition \ref{equivalent} and Theorem
\ref{thm3- spectral measure}.
\end{proof}

\section{Proof of Theorem \ref{thm-main3}}\label{sect.5}

Throughout this section, the symbol \(\nu\) is reused for the general
rational-scale infinite convolution defined below; it does not necessarily
refer to the particular measure introduced in Section~\ref{sect.4}.

Let
\(d_k=l_k/t_k>0\) be written in lowest terms and let \(\gamma_k\geq2\).
Assume that
\begin{equation}\label{rational-convergence}
 \sum_{k=1}^{\infty}
 \frac{\gamma_k-1}{d_1d_2\cdots d_k}<\infty,
\end{equation}
so that
\begin{equation}\label{new moran}
 \nu=\delta_{d_1^{-1}D_{\gamma_1}}*
 \delta_{(d_1d_2)^{-1}D_{\gamma_2}}*\cdots*
 \delta_{(d_1\cdots d_k)^{-1}D_{\gamma_k}}*\cdots
\end{equation}
is a compactly supported probability measure. Write
\[
 \omega_{>k}=\delta_{d_{k+1}^{-1}D_{\gamma_{k+1}}}*
 \delta_{(d_{k+1}d_{k+2})^{-1}D_{\gamma_{k+2}}}*\cdots .
\]

For \(k\geq1\), write
\[
 d_1d_2\cdots d_k=\frac{L_k}{T_k},
 \qquad \gcd(L_k,T_k)=1.
\]
Assume that \(\{T_k\}_{k\geq1}\) and
\(\{\gamma_k\}_{k\geq1}\) are bounded. Choose \(c\in\N\) divisible
by  \(l_1T_k\gamma_k\) for every \(k\). Then the zero set of $\widehat{\nu}$ satisfies 
\begin{equation*}
\mathcal{Z}(\widehat{\nu})=\bigcup_{k=1}^{\infty} d_1\cdots d_k\left(\frac{\mathbb{Z}\setminus \gamma_k\mathbb{Z}}{\gamma_k}\right)\subset\frac{d_1}{c} \mathbb{Z}.
\end{equation*}

Let  $ \Lambda $ be  a spectrum of $\nu$ with $0\in \Lambda$, we have $({c}/{d_1})\Lambda\subset\mathbb{Z} $. It follows that
\begin{equation*}
\frac{1}{d_1}\Lambda=\bigsqcup_{n=0}^{c-1} \left(\frac{n}{c}+\Lambda_n\right),
\end{equation*}
where $ \Lambda_n=\mathbb{Z}\cap\left(\frac{\Lambda}{d_1}-\frac{n}{c}\right) $and $\frac{n}{c}+\Lambda_n=\emptyset$ if $\Lambda_n=\emptyset$. Moreover, the union is disjoint because the residue \(n\) is uniquely determined
modulo \(c\).

Denote by $ q_k:={c}/{\gamma_k}\in\mathbb{N}, \  k\geq1 $. For any $ n\in\{0,1,\dots,c-1\} $, there exists a unique pair of integers $ (i,j)\in\{0,1,\dots,q_1-1\}\times \{0,1,\dots,\gamma_1-1\} $ such that
$$
n=i+q_1 j.
$$
Reindexing the preceding disjoint union therefore gives
\begin{equation}\label{decomposition}
	\frac1{d_1}\Lambda
	=
	\bigsqcup_{i=0}^{q_1-1}
	\bigsqcup_{j=0}^{\gamma_1-1}
	\left(
	\frac{i+q_1j}{c}
	+\Lambda_{i+q_1j}
	\right).
\end{equation}

\begin{Lem}\label{5.2}
Let \(\nu\) be defined by \eqref{new moran} under the assumptions above.
If \(\nu\) is spectral with a spectrum \(\Lambda\) containing \(0\), then
\(\omega_{>1}\) is spectral. Moreover, for any
\(j_i\in\{0,1,\ldots,\gamma_1-1\}\), \(0\leq i<q_1\), the set
\[
 \Gamma=\bigsqcup_{i=0}^{q_1-1}
 \left(\frac{i+q_1j_i}{c}+\Lambda_{i+q_1j_i}\right)
\]
is a spectrum of \(\omega_{>1}\), provided that \(\Gamma\neq\varnothing\).
\end{Lem}

\begin{proof}
	For \(0\leq i<q_1\) and \(0\leq j<\gamma_1\), put
	\[
	\Lambda_{i,j}:=\Lambda_{i+q_1j},
	\qquad
	r_{i,j}:=\frac{i+q_1j}{c}.
	\]
	Fix a choice \(\{j_i\}_{i=0}^{q_1-1}\) for which \(\Gamma	=
	\bigsqcup_{i=0}^{q_1-1}
	\left(r_{i,j_i}+\Lambda_{i,j_i}\right)\)	is nonempty.
	
	We first prove that \(\Gamma\) is a bi-zero set of \(\omega_{>1}\). By the
	definition, we have \(d_1\Gamma\subset\Lambda\).  Let \(\lambda_1,\lambda_2\in\Gamma\) be distinct. Write
	\[
	\lambda_s=r_{i_s,j_{i_s}}+z_s,
	\qquad
	z_s\in\Lambda_{i_s,j_{i_s}},
	\qquad s=1,2.
	\]
  Since
	\(d_1\lambda_1,d_1\lambda_2\in\Lambda\), orthogonality of \(E(\Lambda)\)
	and the identity \(\widehat\nu(d_1t)	= 	m_{D_{\gamma_1}}(t)\widehat\omega_{>1}(t)\) 	 give
	\begin{equation}\label{mask-tail-zero}
		m_{D_{\gamma_1}}(\lambda_1-\lambda_2)\widehat\omega_{>1}(\lambda_1-\lambda_2)=0.
	\end{equation}
	Recall that \(m_{D_{\gamma_1}}\) is \(1\)-periodic,
	\[
	m_{D_{\gamma_1}}(n)=1\quad(n\in\mathbb Z),
	\qquad
	\mathcal Z(m_{D_{\gamma_1}})
	=
	\frac{\mathbb Z\setminus\gamma_1\mathbb Z}{\gamma_1}.
	\]
	
	If \(i_1=i_2\), then \(r_{i_1,j_{i_1}}=r_{i_2,j_{i_2}}\), so
	\(\lambda_1-\lambda_2=z_1-z_2\in\mathbb Z\). Thus the mask factor in
	\eqref{mask-tail-zero} equals \(1\), and
	\(\widehat\omega_{>1}(\lambda_1-\lambda_2)=0\).  Suppose that \(i_1\neq i_2\). Since \(z_1-z_2\in\mathbb Z\), periodicity
	gives
	\[
	m_{D_{\gamma_1}}(\lambda_1-\lambda_2)
	=
	m_{D_{\gamma_1}}
	\left(r_{i_1,j_{i_1}}-r_{i_2,j_{i_2}}\right).
	\]
	If this quantity were zero, then for some \(\ell\in\mathbb Z\),
	\[
	\frac{i_1-i_2+q_1(j_{i_1}-j_{i_2})}{c}
	=
	\frac{\ell}{\gamma_1}.
	\]
	Since \(c=q_1\gamma_1\), it would follow that
	\(i_1-i_2	= 	q_1\bigl(\ell-j_{i_1}+j_{i_2}\bigr)	\in q_1\mathbb Z\),  	contrary to  \(0<|i_1-i_2|<q_1\). 	Thus the mask factor is nonzero, and \eqref{mask-tail-zero} again yields
	\(\widehat\omega_{>1}(\lambda_1-\lambda_2)=0\). Therefore \(\Gamma\) is a bi-zero set
	of \(\omega_{>1}\).

It remains to prove completeness. For \(t\in\R\), define
\[
\begin{aligned}
	p_{i,j}(t)
	&:=
	\left|m_{D_{\gamma_1}}(t+r_{i,j})\right|^2,\\
	x_{i,j}(t)
	&:=
	\sum_{z\in\Lambda_{i,j}}
	\left|
	\widehat\omega_{>1}(t+r_{i,j}+z)
	\right|^2.
\end{aligned}
\]
Each component \(r_{i,j}+\Lambda_{i,j}\) is a bi-zero set of
\(\omega_{>1}\). Hence \(x_{i,j}\) is continuous by the observation
following Theorem \ref{J-P}.

By \eqref{decomposition}, every element of \(d_1^{-1}\Lambda\) has a unique
representation \(r_{i,j}+z\), where \(z\in\Lambda_{i,j}\). Moreover,
\[
\widehat\nu\bigl(d_1(t+r_{i,j}+z)\bigr)
=
m_{D_{\gamma_1}}(t+r_{i,j})
\widehat\omega_{>1}(t+r_{i,j}+z),
\]
because \(z\in\mathbb Z\) and the mask polynomial is \(1\)-periodic.
Summing over the disjoint decomposition of \(d_1^{-1}\Lambda\) and using
\(Q_{\Lambda,\nu}\equiv1\), we obtain
\begin{equation}\label{weighted-tail-identity}
	1
	=
	Q_{\Lambda,\nu}(d_1t)
	=
	\sum_{i=0}^{q_1-1}
	\sum_{j=0}^{\gamma_1-1}
	p_{i,j}(t)x_{i,j}(t).
\end{equation}

For each fixed \(i\), the set \(\{i+q_1j:0\leq j<\gamma_1\}\)
is a spectrum of \(\delta_{c^{-1}D_{\gamma_1}}\). Therefore,
\[
\sum_{j=0}^{\gamma_1-1}p_{i,j}(t)=1.
\]

Suppose now that \(t\notin\mathbb Q\). Since \(r_{i,j}\in\mathbb Q\), all zeros of \(m_{D_{\gamma_1}}\) are rational, and all \(p_{i,j}(t)>0\). Choose, for  each  \(i\), an index \(j\in\{0,1,\ldots,\gamma_1-1\}\) for which \(x_{i,j}(t)\) is maximal,
and select the corresponding component \(r_{i,j}+\Lambda_{i,j}\). By the preceding argument, the union of
the selected components is a bi-zero set of \(\omega_{>1}\). Moreover, the
union is disjoint, since components corresponding to different indices \(i\)
have different fractional parts modulo \(1\). Hence the exponentials
\(e_{-\lambda}\) indexed by this union form an orthonormal family in
\(L^2(\omega_{>1})\). Applying Bessel's inequality to \(e_t\), we obtain
\[
\sum_{i=0}^{q_1-1}
\max_{0\leq j<\gamma_1}x_{i,j}(t)
\leq
\|e_t\|_{L^2(\omega_{>1})}^2
=1.
\]

Applying  Lemma \ref{combined sum} to \eqref{weighted-tail-identity} gives
\begin{equation}\label{5.4.0}
	x_{i,0}(t)=x_{i,1}(t)=\cdots=x_{i,\gamma_1-1}(t)
	\quad(0\leq i<q_1),
	\qquad
	\sum_{i=0}^{q_1-1}x_{i,0}(t)=1.
\end{equation}
 Since all \(x_{i,j}\) are continuous and \(\R\setminus\mathbb Q\) is dense in
\(\R\), \eqref{5.4.0} holds for every \(t\in\R\).

Returning to the fixed choice \(\{j_i\}\) defining \(\Gamma\), we obtain
\[
Q_{\Gamma,\omega_{>1}}(t)
=
\sum_{i=0}^{q_1-1}x_{i,j_i}(t)
=
\sum_{i=0}^{q_1-1}x_{i,0}(t)
=1
\qquad(t\in\R).
\]
Together with the bi-zero property, Theorem \ref{J-P} shows that \(\Gamma\) is a spectrum of \(\omega_{>1}\).

Finally, \(0\in\Lambda\), so \(d_1^{-1}\Lambda\neq\varnothing\). Hence the decomposition \eqref{decomposition} contains a
nonempty component. Choosing \(\Gamma\) to include this component shows that
\(\omega_{>1}\) is spectral.
\end{proof}

\begin{Rem}\label{Remark5.2}
For each fixed \(k\), the reduced denominators of
\(d_{k+1}\cdots d_j=(L_j/T_j)/(L_k/T_k)\) are bounded as \(j\) varies.
Hence Lemma \ref{5.2} applies to every tail. The normalizing integer may depend on
the tail; when finitely many successive decompositions are used, we choose one
common multiple that is admissible for all of them.

More precisely, let \(\Gamma_k\) be a spectrum of \(\omega_{>k}\) containing
\(0\). Choose an admissible integer  \(c_k\) for this tail and put
\(q_{k+1}=c_k/\gamma_{k+1}\). For
\(j_i\in\{0,1,\ldots,\gamma_{k+1}-1\}\), \(0\leq i<q_{k+1}\), the set
\[
 \Gamma_{k+1}=\bigsqcup_{i=0}^{q_{k+1}-1}
 \left(\frac{i+q_{k+1}j_i}{c_k}+\Lambda_{i+q_{k+1}j_i}\right)
\]
is a spectrum of \(\omega_{>k+1}\) whenever it is nonempty, where
\[
 \Lambda_{i+q_{k+1}j_i}=\mathbb Z\cap
 \left(\frac{\Gamma_k}{d_{k+1}}-
 \frac{i+q_{k+1}j_i}{c_k}\right).
\]
\end{Rem}

\begin{Prop}\label{5.3New}
Assume that the measure \(\nu\) in \eqref{new moran} is spectral and satisfies
the hypotheses of Lemma \ref{5.2}. Then, for every \(i\geq1\),
\begin{enumerate}[(i)]
\item if \(\gamma_i\nmid t_{i+1}\), then
\(\gamma_{i+1}\mid l_{i+1}\);
\item if \(\gamma_i\nmid t_{i+1}\), \(d_{i+2}=l_{i+2}\in\N\), and
\(t_{i+1}\mid l_{i+2}\), then
\(t_{i+1}\gamma_{i+2}\mid l_{i+2}\).
\end{enumerate}
\end{Prop}

\begin{proof}
By Remark \ref{Remark5.2}, it is enough to prove the assertions for \(i=1\).
Choose one integer \(c\) admissible for the first three tails and write
\(q_r=c/\gamma_r\), \(r=1,2,3\). Let \(\Lambda\) be a spectrum of
\(\nu\) with \(0\in\Lambda\), and put
\[
 \Lambda_{i,j}=\mathbb Z\cap
 \left(\frac{\Lambda}{d_1}-\frac{i+q_1j}{c}\right).
\]

(i) Assume that \(\gamma_1\nmid t_2\). Set \(j_0=0\). For
\(1\leq i<q_1\), choose
\[
 j_i=
 \begin{cases}
 0,&\displaystyle \frac{l_2}{\gamma_2}
      \not\equiv\frac{t_2i}{c}\pmod1,\\[10pt]
 1,&\displaystyle \frac{l_2}{\gamma_2}
      \equiv\frac{t_2i}{c}\pmod1.
 \end{cases}
\]
Changing \(j_i\) from \(0\) to \(1\) changes
\(t_2(i+q_1j_i)/c\) by \(t_2/\gamma_1\), which is not an integer.
Hence
\begin{equation}\label{avoided-congruence}
 \frac{l_2}{\gamma_2}
 \not\equiv t_2\frac{i+q_1j_i}{c}\pmod1
 \qquad(1\leq i<q_1).
\end{equation}
The set
\[
 \Gamma=\bigsqcup_{i=0}^{q_1-1}
 \left(\frac{i+q_1j_i}{c}+\Lambda_{i,j_i}\right)
\]
contains \(0\), and hence is a spectrum of \(\omega_{>1}\) by Lemma
\ref{5.2}. Its next-stage decomposition is
\[
 \Gamma=d_2\bigsqcup_{a=0}^{q_2-1}\bigsqcup_{j=0}^{\gamma_2-1}
 \left(\frac{a+q_2j}{c}+\Gamma_{a,j}\right),
 \qquad
 \Gamma_{a,j}=\mathbb Z\cap
 \left(\frac{\Gamma}{d_2}-\frac{a+q_2j}{c}\right).
\]
Since \(0\in\Gamma_{0,0}\), the equality of the component sums in
\eqref{5.4.0}, applied to this tail, gives
\(\Gamma_{0,1}\neq\varnothing\). Choose \(z'\in\Gamma_{0,1}\). Comparing
the two displayed decompositions of \(\Gamma\), there exist
\(0\leq i<q_1\) and \(z\in\Lambda_{i,j_i}\) such that
\[
 \frac{l_2}{t_2}\left(\frac1{\gamma_2}+z'\right)
 =\frac{i+q_1j_i}{c}+z.
\]
After multiplication by \(t_2\) and reduction modulo \(1\),
\eqref{avoided-congruence} rules out every \(i\neq0\). Thus \(i=j_0=0\), and
\[
 \frac{l_2}{\gamma_2}+l_2z'=t_2z.
\]
Therefore \(\gamma_2\mid l_2\).

(ii) Assume in addition that \(d_3=l_3\in\N\) and
\(t_2\mid l_3\). Again set \(j_0=0\); for \(i\neq0\), take \(j_i=0\) if
\(t_2i/c\notin\mathbb Z\), and \(j_i=1\) otherwise. Since
\(\gamma_1\nmid t_2\),
\begin{equation}\label{nonintegral-selected-residue}
 t_2\frac{i+q_1j_i}{c}\notin\mathbb Z
 \qquad(1\leq i<q_1).
\end{equation}
Let \(\Gamma\) and \(\Gamma_{a,j}\) be as above for this choice. If
\(z'\in\Gamma_{0,0}\), then \(d_2z'\in\Gamma\), so
\[
 \frac{l_2}{t_2}z'=\frac{i+q_1j_i}{c}+z
\]
for some \(i\) and some \(z\in\Lambda_{i,j_i}\). Multiplication by \(t_2\)
and \eqref{nonintegral-selected-residue} force \(i=j_0=0\). Hence
\(l_2z'=t_2z\). Since \(\gcd(l_2,t_2)=1\), we get \(\Gamma_{0,0}\subset t_2\mathbb Z\).

Choose \(0\leq j'_a<\gamma_2\), with \(j'_0=0\), and set
\[
 \Lambda^*=\bigsqcup_{a=0}^{q_2-1}
 \left(\frac{a+q_2j'_a}{c}+\Gamma_{a,j'_a}\right).
\]
This set contains \(0\) and is therefore a spectrum of \(\omega_{>2}\) by
Lemma \ref{5.2}. At the next stage,
\[
 \Lambda^*=l_3\bigsqcup_{s=0}^{q_3-1}
\bigsqcup_{r=0}^{\gamma_3-1}
 \left(\frac{s+q_3r}{c}+\Lambda^*_{s,r}\right).
\]
Here \(t_3=1\), and therefore \(\gamma_2\nmid t_3\). Part (i), applied at
this stage, gives \(\gamma_3\mid l_3\). Moreover,
\(\Lambda^*_{0,r}\neq\varnothing\) for every \(r\). Choose
\(w_r\in\Lambda^*_{0,r}\) and put
\[
 y_r=l_3\left(\frac r{\gamma_3}+w_r\right).
\]
Then \(y_r\in\Lambda^*\cap\mathbb Z\). In the preceding decomposition of
\(\Lambda^*\), an integer can occur only in the component \(a=0\),
\(j'_0=0\); hence \(y_r\in\Gamma_{0,0}\subset t_2\mathbb Z\).
Writing \(l_3=t_2l_3'\) and taking \(r=1\), we obtain
\[
 l_3'\left(\frac1{\gamma_3}+w_1\right)\in\mathbb Z.
\]
Thus \(\gamma_3\mid l_3'\), which is equivalent to
\(t_2\gamma_3\mid l_3\).
\end{proof}

\begin{Prop}\label{5.1}
Assume that \(2\mid p_k\) for every \(k\geq1\). If \(\nu\) in \eqref{eq-infinite convolution-equivalent} is spectral,
then \(2\mid b_{k+1}\) for every \(k\geq1\).
\end{Prop}

\begin{proof}
Since \(p_k\) is even,
\begin{equation}\label{even-digit-splitting}
 D_{p_k}=D_2\oplus 2D_{p_k/2},
\end{equation}
and hence
\[
 \delta_{B_k^{-1}D_{p_k}}
 =\delta_{B_k^{-1}D_2}*\delta_{2B_k^{-1}D_{p_k/2}}.
\]
Thus the convolution defining \(\nu\) contains the factor
\(\delta_{B_k^{-1}D_2}\) at every level \(k\).

Suppose, to the contrary, that \(b_{n+1}\) is odd for some \(n\geq1\). Then
\begin{equation}\label{nested-binary-zeros}
 \mathcal Z\!\left(\widehat{\delta_{B_{n+1}^{-1}D_2}}\right)=\frac{B_{n+1}}2(2\mathbb Z+1)\subset\frac{B_n}2(2\mathbb Z+1)
 =\mathcal Z\!\left(\widehat{\delta_{B_n^{-1}D_2}}\right),
\end{equation}
because multiplication by the odd integer \(b_{n+1}\) preserves oddness.
Write
\[
 \nu=\nu'*\delta_{B_{n+1}^{-1}D_2},
\]
where \(\nu'\) is the convolution of all remaining factors. By
\eqref{even-digit-splitting},  $\mathcal Z\!\left(\widehat{\delta_{B_n^{-1}D_2}}\right)\subset \mathcal Z(\widehat{\nu'})$, 
 so \eqref{nested-binary-zeros} gives
\[
 \mathcal Z\!\left(\widehat{\delta_{B_{n+1}^{-1}D_2}}\right)
 \subset\mathcal Z(\widehat{\nu'}).
\]
Therefore
\[
 \mathcal Z(\widehat\nu)
 =\mathcal Z(\widehat{\nu'})
  \cup\mathcal Z\!\left(\widehat{\delta_{B_{n+1}^{-1}D_2}}\right)
 =\mathcal Z(\widehat{\nu'}).
\]
If \(\Lambda\) were a spectrum of \(\nu\), every nonzero element of
\(\Lambda-\Lambda\) would therefore be a zero of \(\widehat{\nu'}\); hence
\(\Lambda\) would be a bi-zero set of \(\nu'\). Lemma \ref{Dai}, applied to
\(\nu=\nu'*\delta_{B_{n+1}^{-1}D_2}\), would then imply that \(\Lambda\)
cannot be a spectrum of \(\nu\), a contradiction. Thus every \(b_{k+1}\)
is even.
\end{proof}

To apply the rational-scale decomposition, absorb the coefficients of the
two-point digit sets in \eqref{eq-form of measure nu} into the scale factors.
Define the scale factors by
\begin{equation}\label{b''k}
 \begin{aligned}
 b''_1&=b_1,&
 b''_{2k+1}&=\frac{2}{b_{k+1}(2p_k-1)-1},\\
 b''_2&=b_2,&
 b''_{2k+2}&=\frac{b_{k+2}(b_{k+1}(2p_k-1)-1)}{2},
 \qquad k\geq1,
 \end{aligned}
\end{equation}
and define the corresponding digit cardinalities by
\begin{equation}\label{p''k}
 p''_1=p_1,\qquad p''_2=p_2,\qquad
 p''_{2k+1}=2,\qquad p''_{2k+2}=p_{k+2}.
\end{equation}
Then \(\nu\) in \eqref{eq-infinite convolution-equivalent} can be rewritten as
\begin{equation}\label{nu}
 \nu=\delta_{(b''_1)^{-1}D_{p''_1}}*
 \delta_{(b''_1b''_2)^{-1}D_{p''_2}}*
 \delta_{(b''_1b''_2b''_3)^{-1}D_{p''_3}}*\cdots.
\end{equation}
For this representation, define the normalized tails
\[
 \omega_{>k}=\delta_{(b''_{k+1})^{-1}D_{p''_{k+1}}}*
 \delta_{(b''_{k+1}b''_{k+2})^{-1}D_{p''_{k+2}}}*
 \delta_{(b''_{k+1}b''_{k+2}b''_{k+3})^{-1}D_{p''_{k+3}}}*\cdots.
\]

For \(k\geq1\), direct multiplication gives
\begin{equation}\label{products-b-double-prime}
 \prod_{j=1}^{2k}b_j''=b_1b_2\cdots b_{k+1},\qquad
 \prod_{j=1}^{2k+1}b_j''=
 \frac{2b_1b_2\cdots b_{k+1}}{b_{k+1}(2p_k-1)-1}.
\end{equation}
The even products in \eqref{products-b-double-prime} are integers, while the
reduced denominator of an odd product divides
\(b_{k+1}(2p_k-1)-1\). Hence boundedness of \(\{b_k\}\) and \(\{p_k\}\) implies
boundedness of the cumulative denominators and of the digit cardinalities
\(p_j''\).  Thus the rational-scale convolution \eqref{nu} satisfies all the hypotheses
of Lemma \ref{5.2}.

\begin{theorem}\label{necessity}
Assume that \(2\mid p_k\) for every \(k\geq1\) and that
\(\{b_k\}_{k\geq1}\) is bounded. If \(\nu\) in \eqref{eq-infinite convolution-equivalent} is spectral, then \(p_2\mid b_2\), and \(2p_k\mid b_k\) for \(k\geq3\).
\end{theorem}

\begin{proof}
	Proposition \ref{5.1} gives \(2\mid b_k\) for every \(k\geq2\).  For the first two factors in \eqref{nu}, the notation of Proposition
	\ref{5.3New} gives
	\[
	\gamma_1=p_1,\qquad
	d_2=b_2,\qquad
	t_2=1,\qquad
	\gamma_2=p_2.
	\]
	Since \(p_1\) is even, \(p_1\nmid1\). Proposition \ref{5.3New}(i) therefore
	yields  \(p_2\mid b_2\).
	
	By Lemma \ref{5.2} and Remark \ref{Remark5.2}, every successive normalized
	tail of \eqref{nu} is spectral. Fix \(k\geq3\). The normalized tail
	\(\omega_{>2k-5}\) is spectral and begins with the factor indexed by
	\(2k-4\). Its first three digit cardinalities in \eqref{p''k} are
	\[
	\gamma_1=p''_{2k-4}=p_{k-1},
	\qquad
	\gamma_2=p''_{2k-3}=2,
	\qquad
	\gamma_3=p''_{2k-2}=p_k.
	\]
	The second and third scale factors of this tail in \eqref{b''k} are
	\[
	d_2=b''_{2k-3}=\frac{2}{a_{k-2}},
	\qquad
	d_3=b''_{2k-2}=\frac{b_ka_{k-2}}2,
	\]
	where \(a_{k-2}=b_{k-1}(2p_{k-2}-1)-1\).
	
	Since \(b_{k-1}\) is even, \(a_{k-2}\) is odd. Thus
	\(d_2=2/a_{k-2}\) is in lowest terms with \(l_2=2, \ t_2=a_{k-2}\).
	Moreover, \(\gamma_1=p_{k-1}\) is even, whereas \(t_2\) is odd, and hence
	\(\gamma_1\nmid t_2\). 	Since \(b_k\) is even,  \(d_3=\frac{b_k}{2}a_{k-2}\in\mathbb N\), and \(t_2\mid d_3\).
	Proposition \ref{5.3New}(ii) therefore applies and gives \(t_2\gamma_3\mid d_3\). That is, \(a_{k-2}p_k \mid \frac{b_ka_{k-2}}2\).
	Cancelling \(a_{k-2}\) yields \(p_k\mid b_k/2\), equivalently \(2p_k\mid b_k\). Since \(k\geq3\) was arbitrary, the proof is complete.
\end{proof}

\bigskip

\begin{proof}[Proof of Theorem \ref{thm-main3}]
The result follows from Proposition \ref{equivalent}, Theorem
\ref{thm3- spectral measure}, and Theorem \ref{necessity}.
\end{proof}


\begin{thebibliography}{99}

\bibitem{An-Fu-Lai2019}
L.-X. An, X. Fu, and C.-K. Lai,
\emph{On spectral Cantor--Moran measures and a variant of Bourgain's sum of
sine problem}, Adv. Math. 349 (2019), 84--124.

\bibitem{An-He2014}
L.-X. An and X.-G. He,
\emph{A class of spectral Moran measures},
J. Funct. Anal. 266 (2014), no.~1, 343--354.

\bibitem{An-He-Lau2015}
L.-X. An, X.-G. He, and K.-S. Lau,
\emph{Spectrality of a class of infinite convolutions},
Adv. Math. 283 (2015), 362--376.

\bibitem{An-He-Li2015}
L.-X. An, X.-G. He, and H.-X. Li,
\emph{Spectrality of infinite Bernoulli convolutions},
J. Funct. Anal. 269 (2015), 1571--1590.

\bibitem{An-Li-Zhang2025}
L.-X. An, Q. Li, and M.-M. Zhang,
\emph{The generalized Fuglede's conjecture holds for a class of
Cantor--Moran measures},
Pacific J. Math. 334 (2025), no.~2, 189--209.

\bibitem{Cao-Deng-Li-Wu2024}
Y.-S. Cao, Q.-R. Deng, M.-T. Li, and S. Wu,
\emph{A note on the spectrality of Moran-type Bernoulli convolutions by Deng
and Li},
Bull. Malays. Math. Sci. Soc. 47 (2024), Paper No.~127.

\bibitem{Dai2012}
X.-R. Dai,
\emph{When does a Bernoulli convolution admit a spectrum?},
Adv. Math. 231 (2012), no.~3--4, 1681--1693.

\bibitem{Dai-He-Lau2014}
X.-R. Dai, X.-G. He, and K.-S. Lau,
\emph{On spectral \(N\)-Bernoulli measures},
Adv. Math. 259 (2014), 511--531.

\bibitem{Deng2022}
Q.-R. Deng and M.-T. Li,
\emph{Spectrality of Moran self-similar measures on \(\mathbb R\)},
J. Math. Anal. Appl. 506 (2022), no.~1, Paper No.~125547.

\bibitem{Deng2023}
Q.-R. Deng and M.-T. Li,
\emph{Spectrality of Moran-type Bernoulli convolutions},
Bull. Malays. Math. Sci. Soc. 46 (2023), Paper No.~136.

\bibitem{Dutkay-Han-Lai2019}
D.-E. Dutkay, J. Haussermann, and C.-K. Lai,
\emph{Hadamard triples generate self-affine spectral measures},
Trans. Amer. Math. Soc. 371 (2019), no.~2, 1439--1481.

\bibitem{Fuglede1974}
B. Fuglede,
\emph{Commuting self-adjoint partial differential operators and a group
theoretic problem},
J. Funct. Anal. 16 (1974), 101--121.

\bibitem{He-He2017}
L. He and X.-G. He,
\emph{On the Fourier orthonormal bases of Cantor--Moran measures},
J. Funct. Anal. 272 (2017), no.~5, 1980--2004.

\bibitem{Hu-Lau2008}
T.-Y. Hu and K.-S. Lau,
\emph{Spectral property of the Bernoulli convolutions},
Adv. Math. 219 (2008), no.~2, 554--567.

\bibitem{Jorgensen-Pedersen1998}
P.-E.-T. Jorgensen and S. Pedersen,
\emph{Dense analytic subspaces in fractal \(L^2\)-spaces},
J. Anal. Math. 75 (1998), 185--228.

\bibitem{Laba-Wang2002}
I.~\L aba and Y. Wang,
\emph{On spectral Cantor measures},
J. Funct. Anal. 193 (2002), no.~2, 409--420.

\bibitem{Li-Miao-Wang2022}
W.-X. Li, J.-J. Miao, and Z.-Q. Wang,
\emph{Weak convergence and spectrality of infinite convolutions},
Adv. Math. 404 (2022), Paper No.~108425.

\bibitem{Li-Miao-Wang2024-1}
W.-X. Li, J.-J. Miao, and Z.-Q. Wang,
\emph{Spectrality of random convolutions generated by finitely many Hadamard
triples},
Nonlinearity 37 (2024), no.~1, Paper No.~015003.

\bibitem{Li-Miao-Wang2024-2}
W.-X. Li, J.-J. Miao, and Z.-Q. Wang,
\emph{Spectrality of infinite convolutions and random convolutions},
J. Funct. Anal. 287 (2024), no.~7, Paper No.~110539.

\bibitem{Liu-Wang2025}
J.-C. Liu and J.-J. Wang,
\emph{A class of spectral measures with \(m\)-alternate contraction ratios in
\(\mathbb R\)},
preprint, arXiv:2510.27322 (2025).

\bibitem{Shi2019}
R. Shi,
\emph{Spectrality of a class of Cantor--Moran measures},
J. Funct. Anal. 276 (2019), no.~12, 3767--3794.

\bibitem{Strichartz2000}
R.-S. Strichartz,
\emph{Mock Fourier series and transforms associated with certain Cantor
measures},
J. Anal. Math. 81 (2000), 209--238.

\bibitem{Wu-Liu2022}
H.-H. Wu and J.-C. Liu,
\emph{Construction of a class of spectral measures},
Forum Math. 34 (2022), no.~6, 1507--1517.

\bibitem{Wu2024}
H.-H. Wu,
\emph{Spectral self-similar measures with alternate contraction ratios and
consecutive digits},
Adv. Math. 443 (2024), Paper No.~109585.

\bibitem{Wu-Xiao2024}
S. Wu and Y.-Q. Xiao,
\emph{Spectrality of a class of infinite convolutions on \(\mathbb R\)},
Nonlinearity 37 (2024), no.~5, Paper No.~055015.

\end{thebibliography}
\end{document}